\documentclass[a4paper,11pt,twoside,english]{article}
\usepackage{babel}
\usepackage{latexsym,amsfonts}
\usepackage{amsthm}
\usepackage{amsmath}
\usepackage{mathtools}
\usepackage{paralist}

\usepackage{lastpage}
\usepackage{fancyhdr}
\fancypagestyle{plain}{

  \fancyhead[L,C]{}
  \fancyhead[R]{\today}
  \fancyfoot[R,L]{}
  \fancyfoot[C]{\thepage \ of \pageref{LastPage}}
}

\title{Some convexity properties in direct integrals and K\"othe-Bochner spaces}
\author{Jan-David Hardtke}
\date{}

\providecommand{\sm}{\setminus}

\providecommand{\N}{\ensuremath{\mathbb{N}}}
\providecommand{\R}{\ensuremath{\mathbb{R}}}

\providecommand{\A}{\ensuremath{\mathcal{A}}}

\providecommand{\eps}{\ensuremath{\varepsilon}}

\providecommand{\comment}[1]{
{\let\thefootnote=\relax
\footnote{#1}}
\addtocounter{footnote}{-1}
}

\providecommand{\keywords}[1]{
{\let\thefootnote=\relax
\footnote{{\em Keywords}: #1}}
\addtocounter{footnote}{-1}
}

\providecommand{\AMS}[1]{
{\let\thefootnote=\relax
\footnote{{\em AMS Subject Classification} (2010): #1}}
\addtocounter{footnote}{-1}
}

\providecommand{\address}{
{\sc \noindent Department of Mathematics \\
Freie Universit\"at Berlin \\
Arnimallee 6, 14195 Berlin \\
Germany \\}
}

\DeclarePairedDelimiter{\set}{\lbrace}{\rbrace}

\DeclarePairedDelimiter{\abs}{\lvert}{\rvert}
\DeclarePairedDelimiter{\norm}{\lVert}{\rVert}

\theoremstyle{definition}
\newtheorem{definition}{Definition}[section]
\newtheorem*{definition*}{Definition}

\newtheorem*{example*}{Example}

\newtheorem*{remark*}{Remark}
\theoremstyle{plain}
\newtheorem{lemma}[definition]{Lemma}
\newtheorem*{lemma*}{Lemma}
\newtheorem{proposition}[definition]{Proposition}
\newtheorem*{proposition*}{Proposition}
\newtheorem{theorem}[definition]{Theorem}
\newtheorem*{theorem*}{Theorem}

\newtheorem*{corolary*}{Corollary}

\newenvironment{Proof}[1][\proofname]{\begin{proof}[#1] \setlength{\parindent}{0pt}}{\end{proof}}
\newenvironment{Abstract}{\centering\begin{minipage}{0.8\textwidth} \noindent \small {\sc Abstract.}}{\end{minipage}\par}

\usepackage{color}
\definecolor{darkgreen}{rgb}{0,0.5,0}

\numberwithin{equation}{section}
\newtagform{colored}[\color{blue}]{\color{blue}(}{\color{blue})}
\usetagform{colored}

\usepackage[colorlinks,linkcolor=blue,citecolor=red,urlcolor=darkgreen]{hyperref}
\providecommand{\email}{{\it E-mail address:} \href{mailto:hardtke@math.fu-berlin.de}{\tt hardtke@math.fu-berlin.de}}

\usepackage{amsrefs}

\begin{document}

\maketitle

\begin{Abstract}
The notion of direct integrals introduced in \cite{haydon} is a generalisation of the well-known concept of K\"othe-Bochner spaces of vector-valued functions (using a family of target spaces instead of just one space). Here we will discuss some classical geometric properties like strict convexity, local uniform convexity and uniform convexity in direct integrals. We will also consider strongly convex and very convex K\"othe-Bochner spaces.
\end{Abstract}
\comment{This work is financed by the Deutsche Forschungsgemeinschaft (DFG), grant number HA 8071/1-1.} 
\AMS{46E40 46E30 46B20}

\section{Introduction}\label{sec:intro}
Throughout this paper we will denote by $(S,\A,\mu)$ a complete, $\sigma$-finite measure space
and by $E$ a K\"othe function space over $(S,\A,\mu)$, that is, $E$ is a Banach space of  real-valued measurable functions on $S$ (modulo equality $\mu$-almost everywhere) such that the following conditions are satisfied:
\begin{enumerate}[(i)]
\item $\chi_A\in E$ for every $A\in \A$ with $\mu(A)<\infty$ (where $\chi_A$ denotes the characteristic function of $A$),
\item for every $f\in E$ and every set $A\in \A$ with $\mu(A)<\infty$ $f$ is $\mu$-integrable over $A$,
\item if $g$ is measurable and $f\in E$ such that $\abs*{g(t)}\leq\abs*{f(t)}$ $\mu$-a.\,e. then $g\in E$ and $\norm{g}_E\leq\norm{f}_E$,
\item there exists a function $f\in E$ such that $f(s)>0$ for $\mu$-a.\,e. $s\in S$.
\end{enumerate}
As classical examples one can consider the spaces $L^p(\mu)$ for $1\leq p\leq\infty$.\par
If $X$ is a Banach space, a function $f:S \rightarrow X$ is called simple if there are finitely many pairwise disjoint sets $A_1,\dots ,A_n\in \A$ such that $\mu(A_i)<\infty$ for all $i=1,\dots,n$, $f$ is constant on each $A_i$ and $f(t)=0$ for every 
$t\in S\sm\bigcup_{i=1}^nA_i$.\par 
The function $f$ is called Bochner-measurable if there is a sequence of simple functions
which converges pointwise almost everywhere to $f$ (in the norm of $X$).\par 
By $E(X)$ we denote the space of all $X$-valued Bochner-measurable functions $f$ on $S$
(modulo equality almost everywhere) such that $\norm{f(\cdot)}\in E$. A norm on $E(X)$ is
defined by $\norm{f}_{E(X)}:=\norm{\norm{f(\cdot)}}_E$. $E(X)$ is again a Banach space,
the so called K\"othe-Bochner space induced by $E$ and $X$.\par 
A generalisation of this concept, namely the notion of direct integrals, was introduced by Haydon, Levy and Raynaud in \cite{haydon}. Let $X$ be a real vector space and $(\norm{\cdot}_s)_{s\in S}$ a family of norms on $X$ such that for each $x\in X$ the function $s \mapsto \norm{x}_s$ is measurable. Denote by $(X_s,\norm{\cdot}_s)$ the completion of $(X,\norm{\cdot}_s)$ for every $s\in S$.\par 
Then a function $f\in \prod_{s\in S}X_s$\footnote{That is, $f:S \rightarrow \bigcup_{s\in S}X_s$ with $f(s)\in X_s$ for each $s\in S$.} is called Bochner-measurable if there is a sequence $(f_n)_{n\in \N}$ of $X$-valued simple functions such that $\norm{f_n(s)-f(s)}_s\to 0$ for almost every $s\in S$. Note that if $f$ is Bochner-measurable, then $s \mapsto \norm{f(s)}_s$ is measurable.\par
The direct integral of $(X_s)_{s\in S}$ with respect to $E$ is defined as the space
of all Bochner-measurable functions $f\in \prod_{s\in S}X_s$ such that $(\norm{f(s)}_s)_{s\in S}\in E$, where we again identify two functions in $\prod_{s\in S}X_s$ if they agree almost everywhere. This space is denoted by $(\int_S^{\oplus}X_s\,\text{d}\mu(s))_E$. Endowed with the norm $\norm{f}_{(X_s)}^E:=\norm{(\norm{f(s)}_s)_{s\in S}}_E$ it becomes a Banach space.\par 
Note that for a single Banach space $X$ we have $(\int_S^{\oplus}X\,\text{d}\mu(s))_E=E(X)$.\par
If $\mu$ is the counting measure on $S$, then $(\int_S^{\oplus}X_s\,\text{d}\mu(s))_E$ is just the usual $E$-direct sum $\big[\bigoplus_{s\in S} X_s\big]_E$.\par
Direct integrals of Hilbert spaces and von-Neumann algebras have been known before \cite{haydon} and they have proved very useful in the theory of von-Neumann algebras, see for instance \cite{dixmier}. Also, in \cite{behrends} a slightly different notion of $L^p$-direct integral modules appeared (the main difference to the spaces $(\int_S^{\oplus}X_s\,\text{d}\mu(s))_{L^p}$ in the sense of \cite{haydon} is that 
$S$ carries a topology and the functions $f:S \rightarrow \bigcup_{s\in S}X_s$ are supposed to be such that $s \mapsto \norm{f(s)}_s$ is a {\it continuous} $L^p$-function).\par
In fact, the direct integrals are only a special case of a more general construction of spaces 
$\mathcal{X}_E$, where $\mathcal{X}$ is a so called randomly normed space, but we will not discuss this notion here. Instead, we refer the reader to \cite{haydon} for information on properties of direct integrals/randomly normed spaces and their applications in Banach space theory. For more information on K\"othe-Bochner spaces, one may consult the book \cite{lin} by Lin.\par
Here we would like to discuss some classical convexity properties in direct integrals. Let us briefly recall the definitions. For a Banach space $X$, we denote by $B_X$ its closed unit ball, by $S_X$ its unit sphere and by $X^*$ its dual space. $X$ is called strictly convex (SC) if $\norm{x+y}<2$ for all $x,y\in S_X$ with $x\neq y$.\par 
$X$ is called locally uniformly convex (LUC) if for every sequence $(x_n)_{n\in \N}$ in $S_X$ and every $x\in S_X$ with $\norm{x_n+x}\to 2$ one has $\norm{x_n-x}\to 0$.\par
Finally, $X$ is uniformly convex (UC) if for all sequences $(x_n)_{n\in \N}, (y_n)_{n\in \N}$ in $S_X$ with $\norm{x_n+y_n}\to 2$ one has $\norm{x_n-y_n}\to 0$.\par 
The latter fact can also be expressed in terms of the modulus of convexity of $X$, which is defined by
\begin{equation*}
\delta_X(\eps):=\inf\set*{1-\norm*{\frac{x+y}{2}}:x,y\in B_X, \norm{x-y}\geq \eps} \ \ \forall \eps\in (0,2].
\end{equation*}
$X$ is uniformly convex if and only if $\delta_X(\eps)>0$ for every $\eps\in (0,2]$.\par
Also, $X$ is said to be midpoint locally uniformly convex (MLUC) if the following holds:
whenever $x\in X$ and $(x_n)_{n\in \N}$ is a sequence in $X$ such that $\norm{x_n\pm x}\to \norm{x}$, then $\norm{x_n}\to 0$ (equivalently, whenever $(x_n)_{n\in \N}$ and $(y_n)_{n\in \N}$ are two sequences in $S_X$ and $x\in S_X$ such that $\norm{x_n+y_n-2x}\to 0$, then
$\norm{x_n-y_n}\to 0$).\par
It is known (see for instance Corollary 5 in the work \cite{hudzik} of Hudzik and Wla\'zlak) that a K\"othe-Bochner space $E(X)$ is SC/LUC/MLUC if and only if both $E$ and $X$ are SC/LUC/MLUC. All these results have been known before \cite{hudzik} (see the references therein), but the proof-technique of \cite{hudzik} using sublinear operators gives even more general results, which, as we shall see in Section 2, are also applicable to direct integrals.\par
Concerning uniform convexity, M. Day proved in \cite{day} that for $1<p<\infty$ the $\ell^p$-sum
of a family $(X_s)_{s\in S}$ of Banach spaces is UC if $\inf_{s\in S}\delta_{X_s}(\eps)>0$
for each $\eps$ and that the Lebesgue-Bochner space $L^p(\mu,X)$ is UC if $X$ is UC. In \cite{day2} he generalised his result to arbitrary direct sums with respect to a proper function space (in our language: a K\"othe space $E$ over a set with the counting measure; the $E$-direct sum $(X_s)_{s\in S}$ is UC if $E$ is UC and 
$\inf_{s\in S}\delta_{X_s}(\eps)>0$ for every $\eps$).\par 
Day also noted in \cite{day2} that his argument generalises further to certain spaces of vector-valued functions, which in our language are K\"othe-Bochner spaces: $E(X)$ is UC if $E$ and $X$ are UC.\par
In \cite{greim}, P. Greim already studied uniform convexity (and uniform smoothness), as well
as strict convexity (and smoothness) for the $L^p$-direct integral modules of \cite{behrends}
that we have mentioned above. He also used Day's technique of \cite{day} for uniform convexity. In the next section we will see that Day's result and its proof also directly generalise to the $E$-direct intergrals of \cite{haydon}.\par
We will also consider two less well-known classes of spaces: very convex and strongly convex spaces. The former were introduced by Sullivan in \cite{sullivan}: $X$ is called very convex (or very rotund) if whenever $x\in S_X$, $x^*\in S_{X^*}$ and $x^{**}\in S_{X^{**}}$ are such that $x^*(x)=1=x^{**}(x^*)$ one already has $x^{**}=x$ (under the canonical identification of $X$ with a subspace of $X^{**}$). In \cite{zhang1} it was proved that $X$ is very convex if and only if the following holds: if $(x_n)_{n\in \N}$ is a sequence in $S_X$, $x\in S_X$ and there exists a functional $x^*\in S_{X^*}$ such that $x^*(x_n)\to 1$ and $x^*(x)=1$, then $(x_n)_{n\in \N}$ is weakly convergent to $x$.\par
The concept of strongly convex spaces was introduced in \cite{wu}: $X$ is called strongly convex if it fulfils the above statement but with weak convergence replaced by norm convergence.\par
In the work \cite{zhang2} it was shown that the class of strongly convex spaces coincides with the class of so called almost locally uniformly rotund (ALUR) spaces and the class of very  convex spaces coincides with that of the so called weakly ALUR spaces.\par 
Concerning stability properties, it was proved in \cite{zhang1} that for $p\in (1,\infty)$ the $p$-direct sum of any family of Banach spaces is very convex if and only if each summand is very convex. In \cite{ren} it was shown that a K\"othe-Bochner space $E(X)$ is strongly convex
if $E$ and $X$ are strongly convex and $X^*$ has the Radon-Nikodym property.\par 
In the next section we will see that this result also holds without the assumption on $X^*$,
and also that $E(X)$ is very convex whenever $X$ is very convex and $E$ is strongly convex.
For this we will need a suitable description of the dual of a K\"othe-Bochner space. First
recall that the K\"othe dual $E^{\prime}$ of $E$ is defined as the space of all measurable functions $g:S\rightarrow \R$ (modulo equality a.\,e.) such that 
\begin{equation*}
\norm{g}_{E^{\prime}}:=\sup\set*{\int_S\abs*{fg}\,\mathrm{d}\mu:f\in B_E}<\infty.
\end{equation*}
Then $(E^{\prime},\norm{\cdot}_{E^{\prime}})$ is again a K\"othe function space and the operator $R:E^{\prime}\rightarrow E^*$ defined by 
\begin{equation*}
(Rg)(f)=\int_S fg\,\mathrm{d}\mu \ \ \forall f\in E, \forall g\in E^{\prime}
\end{equation*}
is well-defined, linear and isometric. Moreover, $R$ is surjective if and only if $E$ is order continuous (see \cite{lin}*{p.149}), thus for order continuous $E$ we have $E^*=E^{\prime}$.\footnote{Order continuity of $E$ means that for every sequence $(f_n)_{n\in \N}$ in $E$ with $f_{n+1}\leq f_n$ for every $n\in \N$ and $\inf_{n\in \N}f_n=0$ one has
$\norm{f_n-f}_E\to 0$ (here and in the following, $g\leq h$ means $g(s)\leq h(s)$
for almost every $s\in S$).}\par
Now if $E$ is order continuous and $X^*$ has the Radon-Nikodym property, then the mapping $T:E^{\prime}(X^*) \rightarrow E(X)^*$ given by
\begin{equation*}
T(F)(f):=\int_{S}F(t)(f(t))\,\text{d}\mu(t) \ \ \forall f\in E(X), \forall F\in E^{\prime}(X^*)
\end{equation*}
is an isometric isomorphism. This follows for example from the general representation theory in \cite{bukhvalov}. This description was also used in \cite{ren} to obtain the aforementioned result on strong convexity in K\"othe-Bochner spaces.\par
If $X^*$ does not necessarily have the Radon-Nikodym property, the description of $E(X)^*$ is more involved. First, a function $F:S\rightarrow X^*$ is called weak*-measurable if $F(\cdot)(x)$ is measurable for every $x\in X$. We define an equivalence relation on the set of
all weak*-measurable functions by setting $F\sim G$ if and only if for every $x\in X$ $F(t)(x)=G(t)(x)$ a.\,e. and we denote by $E^{\prime}(X^*,w^*)$ the space of all equivalence classes of weak*-measurable functions $F$ such that there is some $g\in E^{\prime}$ with 
$\norm{F(t)}\leq g(t)$ a.\,e.\par
A norm on $E^{\prime}(X^*,w^*)$ is defined by
\begin{equation*}
\norm{[F]}_{E^{\prime}(X^*,w^*)}:=\inf\set*{\norm{g}_{E^{\prime}}:g\in E^{\prime} \ \text{and} \ \norm{F(t)}\leq g(t) \ \text{a.\,e.}}.
\end{equation*}
Then the following deep theorem holds.
\begin{theorem}[\cite{bukhvalov}, see also \cite{lin}*{Theorem 3.2.4.}]\label{thm:dual E(X)}
If $E$ is order continuous and $X$ is any Banach space, then the map $V:E^{\prime}(X^*,w^*)\rightarrow E(X)^*$ defined by
\begin{equation*}
V([F])(f):=\int_S F(t)(f(t))\,\mathrm{d}\mu(t) \ \ \forall f\in E(X), \forall [F]\in E^{\prime}(X^*,w^*)
\end{equation*}
is an isometric isomorphism. Moreover, every equivalence class $L$ in $E^{\prime}(X^*,w^*)$ has a representative $F$
such that $\norm{F(\cdot)}\in E^{\prime}$ and $\norm{L}_{E^{\prime}(X^*,w^*)}=\norm*{\norm{F(\cdot)}}_{E^{\prime}}$.
\end{theorem}
We will use this result in the proof of Theorem \ref{thm:veryconvex}. Also, we will need the following two results on K\"othe function spaces. The first one was used in \cite{ren} and in \cite{hudzik} (Lemma 2 there, see \cite{akilov}*{Lemma 2 on p.97} for a proof).
\begin{lemma}\label{lemma:subsequences}
If $(f_n)_{n\in \N}$ is a sequence in $E$ and $f\in E$ such that $\norm{f_n-f}_E\to 0$,
then there exists a function $g\in E$ with $g\geq 0$, a subsequence $(f_{n_k})_{k\in \N}$ and a sequence $(\eps_k)_{k\in \N}$ in $(0,\infty)$ which decreases to $0$ such that $\abs{f_{n_k}(t)-f(t)}\leq \eps_kg(t)$ a.\,e. for every $k\in \N$. In particular, $f_{n_k}(t)\to f(t)$ a.\,e.
\end{lemma}
The second one is an abstract Lebesgue theorem for order continuous K\"othe function spaces
(see for instance \cite{lin}*{Theorem 3.1.7.}).
\begin{lemma}\label{lemma:lebesgue}
If $E$ is order continuous, $(f_n)_{n\in \N}$ is a sequence in $E$, $f\in E$
such that $f_n\to f$ pointwise a.\,e. and there exists $g\in E$ with $|f_n|\leq g$ for every $n\in \N$, then $\norm{f_n-f}_E\to 0$.
\end{lemma}

\section{Results and proofs}\label{sec:results}
Henceforth we will assume that $X$ is a real vector space endowed with a family $(\norm{\cdot}_s)_{s\in S}$ of norms such the function $s \mapsto \norm{x}_s$ is 
measurable for every $x\in X$, and $(X_s,\norm{\cdot}_s)$ is the completion of $(X,\norm{\cdot}_s)$.\par 
Concerning the properties SC/LUC/MLUC in direct integrals, the following holds true.
\begin{proposition}\label{prop:SC-LUC}
If $E$ is SC/LUC/MLUC and almost every $X_s$ is SC/LUC/\\MLUC, then $(\int_S^{\oplus}X_s\,\text{d}\mu(s))_E$ is also SC/LUC/MLUC.\par
In particular, $(\int_S^{\oplus}X_s\,\text{d}\mu(s))_{L^p}$ is SC/LUC/MLUC if
almost every $X_s$ is SC/LUC/MLUC and $p\in (1,\infty)$.
\end{proposition}
(Compare with \cite{greim}*{Theorem 4} on strict convexity in the $L^p$-direct integral modules of \cite{behrends}.)\par 
As mentioned above, these statements can be easily deduced from the results of \cite{hudzik},
which we shall now describe. Denote by $L^0(\mu)$ the space of all equivelance classes of real-valued measurable functions. Let $Z$ be a real vector space and $S:Z \rightarrow L^0(\mu)$
be a sublinear operator, i.\,e. $S(x+y)\leq Sx+Sy$ and $S(\lambda x)=|\lambda|Sx$ for all
$x,y\in Z$ and all $\lambda\in \R$. Assume further that $Sx\geq 0$ for every $x\in Z$ and $Sx=0$ only if $x=0$.\par
Then one can consider the space $D_E(S):=\set*{x\in Z:Sx\in E}$, equipped with the norm $\norm{x}_{D_E(S)}:=\norm{Sx}_E$.\par 
If $X$ is a Banach space, $Z$ is the space of all (equivalence classes of) $X$-valued Bochner-measurable functions and $S(f):=\norm{f(\cdot)}$ for $f\in Z$, then $D_E(S)=E(X)$.\par 
Likewise, if $Z$ is the space of all (equivalence classes of) Bochner-measurable functions $f\in \prod_{s\in S}X_s$ and $S(f)(s):=\norm{f(s)}_s$, then $D_E(S):=(\int_S^{\oplus}X_s\,\text{d}\mu(s))_E$.\par
For the general setting of spaces $D_E(S)$, the following notions were introduced in \cite{hudzik}:
\begin{enumerate}[\upshape(i)]
\item $S$ is extreme at a point $x\in Z$ if $y\in Z$ and $S(x\pm y)=Sx$ implies $y=0$.
\item $S$ is strongly extreme at $x\in Z$ if for any sequence $(x_n)_{n\in \N}$ in $Z$ 
one has that $S(x_n\pm x)\to Sx$ a.\,e. implies that $Sx_n\to 0$ a.\,e.\par 
\item $S$ is locally uniformly rotund at $x\in Z$ if the following holds for any sequence $(x_n)_{n\in \N}$ in $Z$: if $S(x+x_n)\to 2Sx$ a.\,e. and $Sx_n\to Sx$ a.\,e., then $S(x-x_n)\to 0$ a.\,e.
\end{enumerate}

Then the following results were proved in \cite{hudzik}:
\begin{enumerate}[\upshape(a)]
\item If $E$ is SC and $S$ is extreme at every $x\in Z$, then $D_E(S)$ is SC (\cite{hudzik}*{Corollary 1}).
\item If $E$ is MLUC and $S$ is strongly extreme at every $x\in Z$, then $D_E(S)$ is MLUC
(\cite{hudzik}*{Corollary 2}).
\item If $E$ is LUC and $S$ is locally uniformly rotund at every $x\in Z$, then $D_E(S)$ is 
LUC (\cite{hudzik}*{Corollary 3}).
\end{enumerate}

These results were then used to obtain the aforementioned results on the properties SC/LUC/MLUC
in K\"othe-Bochner spaces (\cite{hudzik}*{Corollary 5}). In the same way, one can deduce Proposition \ref{prop:SC-LUC}: if almost every $X_s$ is SC/MLUC/LUC, then it is easily checked that the sublinear operator $S$ given by $S(f)(s):=\norm{f(s)}_s$ is extreme/strongly extreme/locally uniformly rotund at every point $f$, so the result follows.\par 
Next we will consider uniform convexity in direct integrals. As mentioned in the introduction,
it is possible to generalise Day's results from \cites{day, day2}. The proof is almost exactly the same as in \cites{day, day2}, but we shall present it here again for the reader's convenience.

\begin{theorem}\label{thm:UC}
Suppose that $E$ is uniformly convex and that there is a null set $N$ such that the function $(0,2]\ni \eps \mapsto \delta(\eps):=\inf_{s\in S\sm N}\delta_{X_s}(\eps)$ is strictly positive. Then $(\int_S^{\oplus}X_s\,\text{d}\mu(s))_E$ is uniformly convex.\par 
In particular, under the above assumption on $(X_s)_{s\in S}$, the Banach space $(\int_S^{\oplus}X_s\,\text{d}\mu(s))_{L^p}$ is uniformly convex for $p\in (1,\infty)$.
\end{theorem}
(Compare with \cite{greim}*{Theorem 1} on uniform convexity in the $L^p$-direct integral modules of \cite{behrends}).

\begin{Proof}
Write $Y:=(\int_S^{\oplus}X_s\,\text{d}\mu(s))_E$ for short and let $\eps\in (0,2)$.
Let $\eta:=\min\set*{1/2,\delta(\eps/4)}$ and $\alpha:=\delta_E(3\eta \eps/4)$.\par 
1) Take $f,g\in S_Y$ such that $\norm{f-g}_Y\geq \eps$ and $\norm{f(s)}_s=\norm{g(s)}_s$
for every $s\in S$. Then $\norm{f+g}_Y\leq 2(1-\alpha)$.\par 	
To see this first put $\beta(s):=\norm{f(s)}_s$ and $\gamma(s):=\norm{f(s)-g(s)}_s$ for each $s\in S$. Then $\gamma(s)\leq 2\beta(s)$ for every $s$.\par
Let $R(s):=\delta(\gamma(s)/\beta(s))$ if $\gamma(s)>0$ and $R(s)=0$ if $\gamma(s)=0$. It is easy to see that 
\begin{equation}\label{eq:1}
\norm{f(s)+g(s)}_s\leq 2(1-R(s))\beta(s) \ \ \forall s\in S\sm N.
\end{equation}
Put $A:=\set*{s\in S:4\gamma(s)>\beta(s)\eps}$ and $B:=S\sm A$. Then
$1=\norm{f}_Y\geq \norm{\beta \chi_B}_E\geq \frac{4}{\eps}\norm{\gamma \chi_B}_E$,
hence $\norm{\gamma \chi_B}_E\leq \eps/4$. This implies
\begin{equation}\label{eq:2}
\norm{\gamma \chi_A}_E\geq \norm{\gamma}_E-\norm{\gamma \chi_B}_E\geq \norm{f-g}_Y-\eps/4
\geq 3\eps/4.
\end{equation}
Let $t:=\beta \chi_B$, $t^{\prime}:=\beta \chi_A$ and $t^{\prime\prime}:=(1-2\eta)t^{\prime}$.
Then $0\leq t+t^{\prime\prime}\leq t+t^{\prime}=\beta$ (hence $t+t^{\prime},t+t^{\prime\prime}\in B_E$) and $\norm{t+t^{\prime}-(t+t^{\prime\prime})}_E=
\norm{t^{\prime}-t^{\prime\prime}}_E=2\eta\norm{t^{\prime}}_E\geq \eta \norm{\gamma \chi_A}_E$.\par
Thus it follows from \eqref{eq:2} that $2-\norm{2t+t^{\prime}+t^{\prime\prime}}_E\geq 2\delta_E(3\eta\eps/4)$. Hence 
\begin{equation}\label{eq:3}
\norm{(1-\eta)t^{\prime}+t}_E\leq 1-\delta_E(3\eta\eps/4)=1-\alpha.
\end{equation}
Using \eqref{eq:1} we obtain for $s\in S\sm N$ 
\begin{align*}
&\frac{1}{2}\norm{f(s)+g(s)}_s\leq (1-R(s))\beta(s)\leq (1-R(s))t^{\prime}(s)+t(s) \\
&\leq (1-\eta)t^{\prime}(s)+t(s),
\end{align*}
where it was used that $R(s)=\delta(\gamma(s)/\beta(s))\geq \delta(\eps/4)\geq \eta$
for each $s\in A$.\par
Since $\mu(N)=0$ it follows that $\norm{f+g}_Y\leq 2\norm{(1-\eta)t^{\prime}+t}_E$ and 
thus \eqref{eq:3} implies $\norm{f+g}_Y\leq 2(1-\alpha)$.\par
2) Now choose $0<\omega<\min\set*{\eps,2\alpha}$ and $\tau>0$ such that
$2(1-\tau)>2(1-\alpha)+\omega$ and $\tau\leq \delta_E(\omega)$.\par
Let $f,g\in S_Y$ such that $\norm{f+g}_Y>2(1-\tau)$. To complete the proof it suffices to 
show that $\norm{f-g}_Y\leq 2\eps$.\par
First note that $2(1-\delta_E(\omega))<\norm{f+g}_Y\leq \norm{(\norm{f(s)}_s+\norm{g(s)}_s)_{s\in S}}_E$, hence 
\begin{equation}\label{eq:4}
\norm{(\norm{f(s)}_s-\norm{g(s)}_s)_{s\in S}}_E<\omega.
\end{equation}
Let $h(s):=\frac{\norm{f(s)}_s}{\norm{g(s)}_s}g(s)$ if $g(s)\neq 0$ and $h(s):=f(s)$ if $g(s)=0$. It is easy to see that $h$ is Bochner-measurable and $\norm{h(s)}_s=\norm{f(s)}_s$
for each $s\in S$, thus $h\in S_Y$.\par 
Furthermore, $\norm{h(s)-g(s)}_s=|\norm{f(s)}_s-\norm{g(s)}_s|$ for every $s\in S$, hence 
$\norm{h-g}_Y=\norm{(\norm{f(s)}_s-\norm{g(s)}_s)_{s\in S}}_E<\omega$, by \eqref{eq:4}.\par 
It follows that $\norm{f+h}_Y\geq \norm{f+g}_Y-\norm{g-h}_Y>2(1-\tau)-\omega>2(1-\alpha)$.
By part 1) this implies $\norm{f-h}_Y<\eps$.\par
Since $\norm{h-g}_Y<\omega<\eps$ it follows that $\norm{f-g}_Y<2\eps$.
\end{Proof}

As mentioned in the introduction, $E$-direct integrals are only a special case of a more
general class of spaces $\mathcal{X}_E$, where $\mathcal{X}$ is a so called randomly normed space (see \cite{haydon}). In \cites{guo1, guo2}, Guo and Zeng studied randomly normed spaces
(or even randomly normed modules) and introduced the notions of random strict convexity and
random uniform convexity for them. They proved (in our notation) that $\mathcal{X}_{L^p}$ is strictly/uniformly convex if $\mathcal{X}$ is randomly strictly convex/randomly uniformly convex, where $p\in (1,\infty)$.\par 
It is likely possible to derive the statement of Proposition \ref{prop:SC-LUC} on
strict convexity of $L^p$-direct integrals and the result of Theorem \ref{thm:UC} on
uniform convexity of $L^p$-direct integrals from the results of Guo and Zeng. But the author
found it is easier to work with the proof-techniques used above to obtain the more general results without ever refering to the notion of randomly normed spaces.\par 
Next we turn to the announced result on strongly convex/very convex spaces.\par
\begin{theorem}\label{thm:veryconvex}
Suppose that $E$ is strongly convex. Then $E(X)$ is strongly convex/very convex whenever
$X$ is strongly convex/very convex.
\end{theorem}

\begin{Proof}
We will first prove the statement for very convex spaces and then indicate the necessary
changes for the case of strong convexity. The proof is similar to the one from \cite{ren} for strong convexity under the additional assumption that $X^*$ has the Radon-Nikodym property.\par
First, it must be noted that the strong convexity of $E$ implies its order continuity
(see \cite{ren}).\par
Now assume that $X$ is very convex and take a sequence $(f_n)_{n\in \N}$ in $S_{E(X)}$,
an element $f\in S_{E(X)}$ and a functional $\varphi\in S_{E(X)^*}$ with $\varphi(f)=1$
and $\varphi(f_n)\to 1$. We want to show that $(f_n)_{n\in \N}$ converges weakly to $f$.\par 
Because of Theorem \ref{thm:dual E(X)} we can find an element $L\in E^{\prime}(X^*,w^*)$
such that $\varphi=V(L)$ and a representative $F\in L$ such that $\norm{F(\cdot)}\in E^{\prime}$ and $\norm{\norm{F(\cdot)}}_{E^{\prime}}=\norm{\varphi}=1$.\par
We have
\begin{equation*}
\varphi(f)=\int_S F(t)(f(t))\,\text{d}\mu(t)\leq \int_S \norm{F(t)}\norm{f(t)}\,\text{d}\mu(t)
\leq \norm{\norm{F(\cdot)}}_{E^{\prime}}\norm{f}_{E(X)}=1
\end{equation*}
and thus
\begin{equation}\label{eq:5}
\int_S \norm{F(t)}\norm{f(t)}\,\text{d}\mu(t)=1
\end{equation}
and 
\begin{equation}\label{eq:6}
F(t)(f(t))=\norm{F(t)}\norm{f(t)} \ \ \text{a.\,e.}
\end{equation}
Analogously one can show that
\begin{equation}\label{eq:7}
\lim_{n\to \infty}\int_S \norm{F(t)}\norm{f_n(t)}\,\text{d}\mu(t)=1
\end{equation}
and 
\begin{equation*}
\lim_{n\to \infty}\int_S (\norm{F(t)}\norm{f_n(t)}-F(t)(f_n(t)))\,\text{d}\mu(t)=0.
\end{equation*}
So by passing to a subsequence we may assume that 
\begin{equation}\label{eq:8}
\lim_{n\to \infty}(\norm{F(t)}\norm{f_n(t)}-F(t)(f_n(t)))=0\ \ \text{a.\,e.}
\end{equation}
Since $E$ is strongly convex it follows from \eqref{eq:5} and \eqref{eq:7} that 
$\norm{\norm{f_n(\cdot)}-\norm{f(\cdot)}}_E\to 0$.\par
Because of Lemma \ref{lemma:subsequences} we can, by passing to a further subsequence, assume that
\begin{equation}\label{eq:9}
\lim_{n\to \infty}\norm{f_n(t)}=\norm{f(t)} \ \ \text{a.\,e.}
\end{equation}
and moreover that there is some $g\in E$ with $|\norm{f_n(t)}-\norm{f(t)}|\leq g(t)$
a.\,e. (for every $n$).\par 
From \eqref{eq:8} and \eqref{eq:9} we obtain
\begin{equation}\label{eq:10}
\lim_{n\to \infty}F(t)(f_n(t))=\norm{F(t)}\norm{f(t)} \ \ \text{a.\,e.}
\end{equation}
Now let $A:=\set*{t\in S:F(t)=0}$ and $B:=\set*{t\in A:f(t)\neq 0}$. By \eqref{eq:5} we have 
\begin{equation*}
1=\int_{S\sm A} \norm{F(t)}\norm{f(t)}\,\text{d}\mu(t)\leq \norm{\norm{f(\cdot)}\chi_{S\sm A}}_E
\leq \norm{f}_{E(X)}=1.
\end{equation*}
Hence $\norm{\norm{f(\cdot)}\chi_{S\sm A}}_E=1$ and $2\geq \norm{\norm{f(\cdot)}+\norm{f(\cdot)}\chi_{S\sm A}}_E\geq 2\norm{\norm{f(\cdot)}\chi_{S\sm A}}_E=2$. Since $E$ is in particular strictly convex this implies $\mu(B)=0$.\par 
Because of \eqref{eq:6}, \eqref{eq:9} and \eqref{eq:10}, the assumption that $X$ is very convex implies that $f_n(t)\to f(t)$ weakly for a.\,e. $t\in S\sm A$.\par 
Since $\mu(B)=0$ we obtain that $f_n(t)\to f(t)$ weakly for a.\,e. $t\in S$.\par 
Now take an arbitrary $\psi\in E(X)^*$ and write $\psi=V([G])$ with $\norm{G(\cdot)}\in E^{\prime}$ and $\norm{\norm{G(\cdot)}}_{E^{\prime}}=\norm{\psi}$ (Theorem \ref{thm:dual E(X)}).\par
Then we have $G(t)(f_n(t))\to G(t)(f(t))$ for a.\,e. $t\in S$ and also
$|G(t)(f_n(t))|\leq \norm{G(t)}\norm{f_n(t)}\leq \norm{G(t)}(\norm{f(t)}+g(t))$ a.\,e.\par 
Thus Lebesgue's Theorem implies $\psi(f_n)\to \psi(f)$ and the proof is finished.\par 
If $X$ is even strongly convex one can proceed analogously to obtain $\norm{f_n(t)-f(t)}\to 0$
a.\,e. Since $\norm{f_n(t)-f(t)}\leq \norm{f_n(t)}+\norm{f(t)}\leq g(t)+2\norm{f(t)}$ a.\,e.,
Lemma \ref{lemma:lebesgue} implies $\norm{f_n-f}_{E(X)}\to 0$. So in this case $E(X)$ is even strongly convex.
\end{Proof}

{\large \sc Remark:} The author does not know whether $E(X)$ is very convex if $X$ and $E$ are both merely very convex, nor if the above Theorem can be generalised to direct integrals
(to do so one would need a suitable description of the dual of a direct integral, which up to the author's knowledge is not known so far).\par
\ \\
{\large \sc Acknowledgements:} The author is very grateful to the DFG (Deutsche Forschungsgemeinschaft (German Research Foundation)) for providing the funding for this work
(grant number HA 8071/1-1).

\address
\email


\begin{bibdiv} 
\begin{biblist}
	
\bib{akilov}{book}{
	title={Functional analysis},
	author={Akilov, G. P.},
	author={Kantorovich, L. V.},
	publisher={Pergamon Press},
	edition={2},
	address={Oxford},
	date={1982}
}

\bib{behrends}{book}{
	title={$L^p$-structure in real Banach spaces},
	author={Behrends, E.},
	author={Danckwerts, R.},
	author={Evans, R.},
	author={G\"obel, S.},
	author={Greim, P.},
	author={Meyfarth, K.},
	author={M\"uller, W.},
	series={Lectures Notes Math.},
	volume={613},
	publisher={Springer},
	address={Berlin-Heidelberg},
	date={1977}
}
	
\bib{bukhvalov}{article}{
	title={On an analytic representation of operators with abstract norm},
	author={Bukhvalov, A. V.},
	journal={Soviet Math. Doklady},
	volume={14},
	date={1973},
	pages={197--201}
}
	
\bib{day}{article}{
	title={Some more uniformly convex spaces},
	author={Day, M. M.},
	journal={Bull. Amer. Math. Soc.},
	volume={47},
	number={6},
	date={1941},
	pages={504--507}
}

\bib{day2}{article}{
	title={Uniform convexity III},
	author={Day, M. M.},
	journal={Bull. Amer. Math. Soc.},
	volume={49},
	number={10},
	date={1943},
	pages={745--750}
}

\bib{dixmier}{book}{
	title={Von Neumann algebras},
	author={Dixmier, J.},
	series={North-Holland Mathematical Library},
	volume={27},
	publisher={North-Holland},
	address={Amsterdam},
	date={1981}
}

\bib{greim}{article}{
	title={Some geometric properties of integral modules},
	author={Greim , P.},
	journal={Manuscripta Math.},
	volume={22},
	number={4},
	date={1977},
	pages={311--323}
}

\bib{guo1}{article}{
	title={Random strict convexity and random uniform convexity in random normed modules},
	author={Guo, T. X.},
	author={Zeng, X. L.},
	journal={Nonlinear Anal: Theory, Methods Appl.},
	volume={73},
	number={5},
	date={2010},
	pages={1239--1263}
}

\bib{guo2}{article}{
	title={An $L^0(\mathcal{F},\R)$-valued function’s intermediate value theorem and its applications to random uniform convexity},
	author={Guo, T. X.},
	author={Zeng, X. L.},
	journal={Acta Math. Sinica, Engl. Ser.},
	volume={28},
	number={5},
	date={2012},
	pages={909--924}
}
	
\bib{haydon}{book}{
	title={Randomly normed spaces},
	author={Haydon, R.},
	author={Levy, M.},
	author={Raynaud, Y.},
	series={Travaux en Cours},
	volume={41},
	publisher={Hermann, \'Editeurs des Sciences et des Arts},
	address={Paris},
	date={1991}
}

\bib{hudzik}{article}{
	title={Rotundity properties in Banach spaces via sublinear operators},
	author={Hudzik, H.},
	author={Wla\'zlak, K.},
	journal={Nonlinear Anal.},
	volume={64},
	date={2006},
	pages={1171--1188}
}

\bib{lin}{book}{
	title={K\"othe-Bochner function spaces},
	author={Lin, P. K.},
	publisher={Birkh\"auser},
	address={Boston-Basel-Berlin},
	date={2004}
}

\bib{ren}{article}{
	title={Rotundity in K\"othe-Bochner spaces},
	author={Ren, L. W.},
	author={Feng, G. C.},
	author={Wu, C. X.},
	journal={Adv. Math. (China)},
	volume={35},
	number={3},
	date={2006},
	pages={350--360}
}

\bib{sullivan}{article}{
	title={Geometrical properties determined by the higher duals of a Banach space},
	author={Sullivan, F.},
	journal={Illinois J. Math.},
	volume={21},
	number={2},
	date={1977},
	pages={315--331}
}

\bib{wu}{article}{
	title={Strong convexity in Banach spaces},
	author={Wu, C. X.},
	author={Li, Y. J.},
	journal={Chin. J. Math.},
	volume={13},
	number={1},
	date={1993},
	pages={105--108}
}

\bib{zhang1}{article}{
	title={On very rotund Banach space},
	author={Zhang, Z. H.},
	author={Zhang, C. J.},
	journal={Appl. Math. Mech.},
	volume={21},
	number={8},
	date={2000},
	pages={965--970}
}

\bib{zhang2}{article}{
	title={Some generalizations of locally and weakly locally uniformly convex space},
	author={Zhang, Z. H.},
	author={Liu, C. Y.},
	journal={Nonlinear Anal.},
	volume={74},
	date={2011},
	pages={3896--3902},
}


\end{biblist}
\end{bibdiv}
\end{document}